# Powering the Future of AI: Navigating the Trade-offs for Europe's Energy Transition and Net-Zero Goals

**Mohammad Hemmati [1], Gbemi Oluleye [2], and Vassilis M. Charitopoulos [1]**

[1] Department of Chemical Engineering, Sargent Centre for Process Systems Engineering, University College London (UCL), Torrington Place, WC1E 7JE, London, United Kingdom

[2] Centre for Environmental Policy, Imperial College London, London SW7 1NE, UK

## Abstract

The rapid expansion of AI globally has led to the proliferation of energy-intensive hyperscale data centres (DCs), making them as a structurally challenging component in power system planning and operation. Using a spatially explicit optimisation model of Europe across 21 AI growth scenarios, we systematically quantify additional demand, capacity requirements, emissions, and operational impacts of DCs. Results indicate that AI could drive 73–723 TWh of extra demand by 2050, risking cumulative emissions overshoots of 67–181 $MtCO_2$ between 2030 and 2050. Our analysis indicates that after 2030, the geography of AI infrastructure will be shaped more by firm power and system flexibility than by the mere abundance of clean energy. In moderate scenarios, AI requires an additional of 200 hours of firm generation, which increases LCOE by €35/MWh in key hubs. We show that even under the pessimistic scenarios, existing infrastructure would require +70 GW additional capacity, while under managed growth pathways, this expansion could reach +226 GW. We further find DCs workload dynamics strongly shape energy dispatch, system flexibility, and emissions, while improved efficiency significantly reduces capacity needs, and system peaks. While our findings suggest that net-zero targets for 2050 may be achieved, critical emission risks may appear in the intermediate years, and the EU may compromise its carbon-neutral goals unless policies adapt to this accelerating digital transformation.

With the global average temperature reaching 1.55°C above pre-industrial levels in 2024, overshooting the Paris Agreement's critical limit of 1.5°C, the pursuit to halt global warming becomes urgent rendering the resilient decarbonisation of energy systems crucial[1,2]. In 2019, Europe issued its Green Deal, which sets out an ambitious roadmap to carbon neutrality by 2050 with intermediate targets including a 55% reduction of greenhouse gas (GHG) emissions by 2030[3]. In the race towards net-zero, key factors including technology mix, pace of decarbonisation and demand management have shown to influence heavily the end cost of the forthcoming power systems transformation[4,5]. Focusing on the latter, the recent widespread adoption of Generative AI (GenAI) and cloud computing services that facilitate ongoing digitalisation efforts, has introduced a new significant demand component that should be accounted for in decarbonisation pathways owing to unprecedented growth in computational workloads[6-8]. Widely adopted AI models consume tens of gigawatt-hours for training and continuous inference[9,10]. For instance, a 100-word ChatGPT-4 response uses roughly 0.14 kWh of electricity[11], and with about 285 million daily users generating responses at that scale, the platform's inference operations can be roughly 0.39 TWh per day (~14.4 TWh annually)[12]. In addition, the training of the GPT-4 model is estimated to have consumed around 50 GWh of electricity, similar to the output of a French nuclear reactor running at full capacity for several days [13].

The global digital landscape is undergoing a massive structural shift, with data centres (DCs) now accounting for approximately 1.5% of global electricity demand[14]. In Europe, this sector currently consumes roughly 70 TWh of electricity, but the trajectory is steepening[15]. While the International Energy Agency (IEA) forecasts a rise to 115 TWh by 2030[15,16], other industry leaders suggest a much steeper climb; McKinsey and Independent Commodity Intelligence Services (ICIS) project totals reaching between 150-168 TWh[17-19], with more aggressive estimates suggesting annual demand could soar as high as 287 TWh[20,21].

Beyond these midterm forecasts, the sheer speed of this acceleration is already creating tangible challenges for power grid that need to reconcile the sector operational intensity with existing infrastructure. In Northern Virginia, the largest DC market in the world, planned gas capacity remained stable between 2021 and 2024, but surged by 20% in 2025 to satisfy skyrocketing demand, forcing some utilities to extend the operation of coal-fired plants to maintain grid stability[22,23]. If these trends are mirrored in Europe, the energy-intensive expansion of AI infrastructure could directly undermine its decarbonisation efforts. However, the challenge must also be viewed through an inverse lens regarding the risk of stranded assets. Should AI deployment or its associated electricity demand fall short of current ever-increasing projections significant over-investment in grid upgrades could lead to significant financial inefficiencies and higher electricity prices for end-consumers [24,25].

DC's power demand in Europe remains insufficiently represented in long-term planning frameworks at both national and continental levels. The primary reference for the future European electricity system ENTSO-E Ten-Year Network Development Plan (TYNDP)[26], and widely-used energy systems models such as PyPSA[27], TIMES[28], Calliope[29] and high-RES[30], do not explicitly incorporate the dynamic and potentially structural electricity loads of AI, and transition pathways have been developed without adequately capturing the scale, temporal profile, and locational implications of DCs. Hence, understanding DCs interactions with the European grid is crucial for ensuring reliable electricity supply and sustainable AI development under net-zero targets[31]. To date researchers have explored the role of DCs as grid responsive assets[32] or assessed their energy and climate impacts[33] while others have analysed their potential to provide operational flexibility to power systems[34]. Some studies have even developed frameworks to quantify the multi energy consumption and industry specific interactions of AI[35]. While these previous works provide valuable insights, they do not fully capture the spatial and temporal dynamics of AI driven electricity demand across the European continent.

Building on this existing literature our study advances the AI-energy nexus by using a spatially explicit optimisation model covering 33 European countries to examine how large-scale AI deployment affects planned grid infrastructure planning and net-zero strategies. The model employs a mixed integer linear programming (MILP) formulation to jointly optimise long term capacity expansion and hourly system dispatch. This approach incorporates existing development plans along with operational constraints and climate policy targets. Using historical data, we derive hourly load profiles of hyperscale DCs and optimise the spatial and temporal evolution of AI related capacity across 21 scenarios ranging from pessimistic to optimistic growth through 2050.

This work is structured around four central questions: **1.** Are existing European power development plans adequate for DC demand growth? **2.** Which infrastructure technologies and energy strategies are essential to reconcile rapid AI and what factors enable or constrain DCs deployment across regions? **3.** What if AI growth stalls or collapses after 2030, and what are the structural risks to power grid expansion? **4.** How will AI growth and efficiency gains reshape European generation patterns, cross-border flows, electricity price and emission pathways?

By addressing these questions, this paper provides the first integrated, spatially explicit, and scenario-based analysis of how large-scale AI-driven demand may interact with Europe's decarbonisation plans, offering new insights for policymakers, system operators, and investors.

# Results

## The Baseline case definition

The baseline case (also referred to as the **Without AI**) represents the projected evolution of the European power system without the additional electricity demand from DCs. To ensure the reliability of our modelling framework, this scenario has been rigorously validated and benchmarked against the TYNDP 2024 Global Ambition [23]. The results demonstrate a high degree of alignment with official projections, showing a total installed capacity of 4664 GW (3.7% deviation from TYNDP projection) and 103 GW of additional cross-border transmission capacity (<1% deviation from TYNDP projection) by 2050 to reach a total of 233 GW. This close convergence with pan-European benchmarks confirms the model's robustness as a stable foundation for evaluating the marginal impacts of AI load growth. A comprehensive description of the baseline assumptions, data sources, and detailed model comparisons is provided in the Supplementary Material, Section S.4.

## AI-Driven demand growth and power system transformation

To systematically explore the uncertainties surrounding the energy footprint of AI, we develop 21 distinct growth scenarios. These scenarios are grounded in projections from the IEA[15], ICIS y[18], and McKinsey[19], and ICIS, with seven scenarios defined for each: Base, Lift-off, Headwinds (Head), High-efficiency (High-eff), and three pessimistic scenarios: Deflation, Medium, and End. This multi-scenario framework allows us to first visualise the energy-impact, infrastructure requirements, and emission trajectories, while answering questions about where and how much demand manifests across Europe for all 21 scenarios. To illustrate the workload dynamics of hyperscale data centers and their operational flexibility, we also run the optimization model using two distinct load profiles derived from historical data, each exhibiting completely different patterns in terms of fluctuations and load factor, referred to as the UK profile and the US profile. The detailed assumption parameters for building all scenarios, are documented in the Methods section and Supplementary Material, Section S.2.

Our analysis indicates that even under the most pessimistic AI growth scenario (73 TWh in 2050), Europe's planned generation and storage capacity may be insufficient to meet AI-driven electricity demand. In most scenarios, wind power emerges as the preferred expansion technology by 2030, with an additional 2–50 GW of installed capacity (around 20% offshore and 80% onshore; Supplementary Figures S.7 and S.8), resulting in approximately +80 TWh of additional generation from wind power (Fig. 1a). Combined-cycle gas turbines (CCGT) become a key balancing technology due to their fast-ramping capability. In medium- and high-growth scenarios, the system requires more than 50 TWh additional CCGT generation.

Notably coal and oil generation increase above baseline levels by 2030. In lift-off scenarios, Oil & Coal generation rises by up to +18 TWh, directly contributing to short-term carbon-budget overshooting and delaying decarbonisation progress (Fig. 1.b). This behaviour persists in 17 of the 21 scenarios until approximately 2045. By contrast, in lower-growth scenarios (e.g., High Eff and Head), clean energy resources absorb the additional load. These scenarios show higher generation from wind and solar, alongside roughly +2% additional nuclear output by 2030, reducing the need for gas, and particularly coal. As a result, emissions decline by up to 12 MtCO2 in some scenarios (blue trajectories in Fig. 1.b).

Beyond 2035, the system's need for low-carbon firm resources increases substantially. The system exhibits 4–30 TWh more CCGT with carbon capture (CCGT-CCS) generation during 2035–2045. By 2050, system requires about +4 GW of $H_2$-CCGT, 2–12 GW of nuclear, +20–60 GW of CCGT, and +3 GW of CCGT-CCS (Fig. 1.a, Supplementary Figures S.7 and S.8). These are driven not only by energy needs but also by reserve requirements, as DCs operations results in higher system peaks. As a result, some gas assets remain online through 2050 to satisfy the reserve requirements.

Renewables and battery storage exhibit non-linear pattern. They expand with AI demand in moderate-to-high growth and high-efficiency scenarios but decline toward the late 2040s in most pathways, except the highest-growth scenarios. Consequently, 17 scenarios deliver 80–160 TWh less wind generation than the baseline by 2050. In high-growth scenarios wind onshore remains a key reliable technology, requiring up to 124 GW of additional capacity by 2050. Low solar costs combined with battery expansion are central in the baseline and lift-off pathways (+224 GW solar, +70 GW battery). Pessimistic scenarios show weaker renewable uptake, as earlier investment in firm technologies ($H_2$-CCGT, nuclear) lowers storage demand, resulting in 10–120 GW less battery capacity.

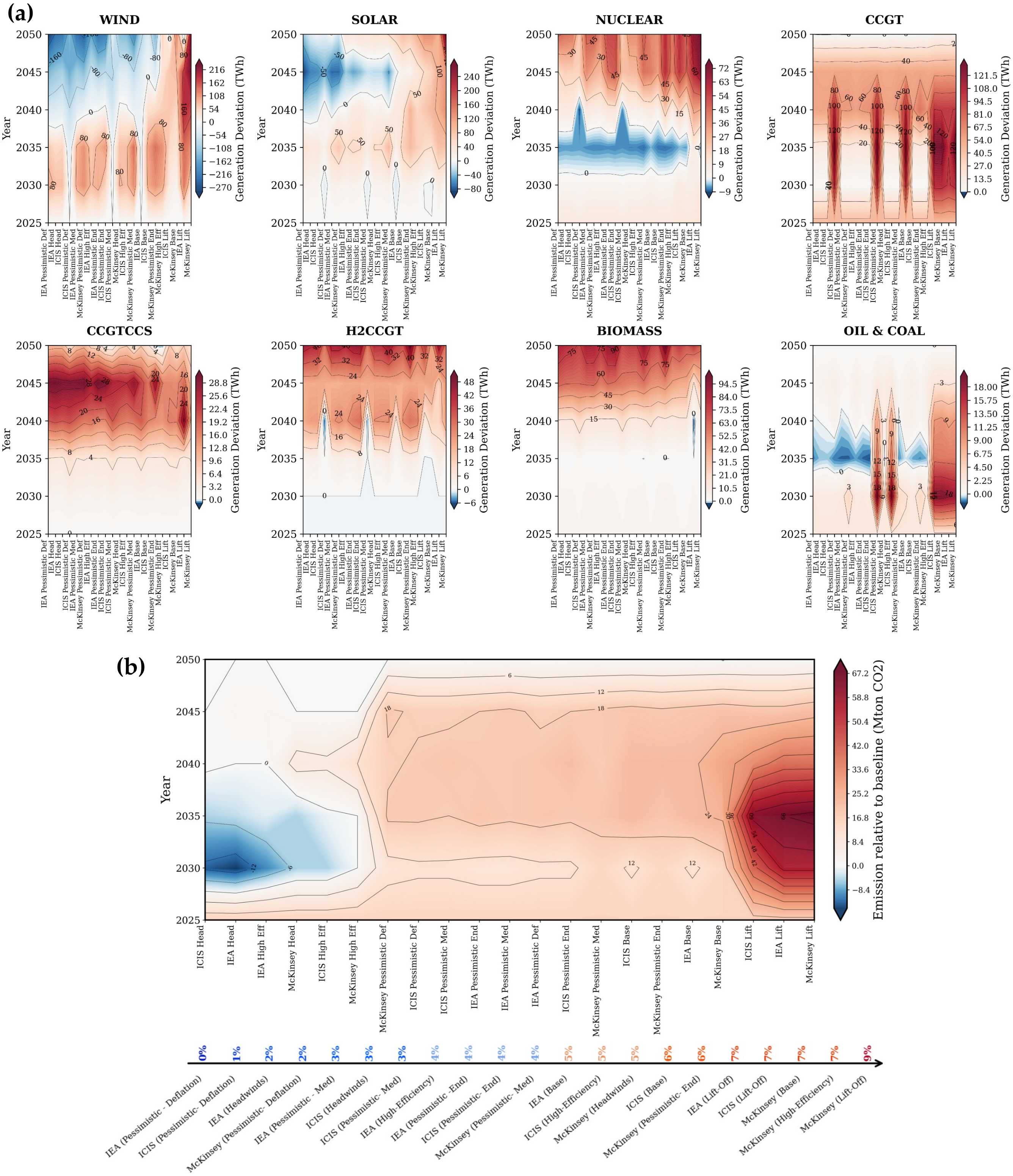

**Fig. 1. Structural evolution of the European power system and emission trajectories under 21 AI-demand scenarios with the US Profile. a,** Heatmaps illustrating the temporal deviation in annual electricity generation (TWh) relative to the baseline scenario (2025–2050) across eight key technologies: Wind, Solar, Nuclear, CCGT, CCGT-CCS, $H_2$-CCGT, Biomass, and Oil & Coal (results on the UK profile is reported in Supplementary Material Figure S.9). The vertical axis represents the transition horizon (Years), while the horizontal axis displays 21 distinct demand scenarios, ordered by their projected CAGR for AI data centers (0% to 9%, indicated by the bottom-most axis). Colour gradients denote generation surplus (red) or deficit (blue) compared to the baseline. **b,** Temporal profile of carbon emission deviations (Mton $CO_2$) across the 21 scenarios. The contour plot highlights the intensity of carbon-budget overshooting (red) or undershooting (blue).

**The evolution of Europe's digital geography: Flexibility as the frontier (post-2030)**

Our analysis indicates that after 2030, the geography of DC infrastructure in Europe will be shaped more by firm power and system flexibility than by the mere abundance of clean energy. Traditional hubs like Germany, the Netherlands, and Ireland are approaching a structural deadlock (Together, these countries host a capacity of 7.5-11.5 GW of DC). During winter peaks, these systems rely heavily on gas-fired generation, incurring severe carbon penalties. The addition of DC exacerbates gas dispatch, leading to a consistent overshooting of emission targets. Consequently, DC growth in these regions plateau post-2030. In the event of collapse scenarios, Germany emerges as the primary candidate for capacity reduction, as the system prefers to curtail DC demand rather than sustain the punitive economic and environmental costs.

Conversely, Portugal and Spain, characterized by mild baseloads and even negative net-load growth[26], offer an ideal location for DC integration. This region leverages significant solar and onshore wind potential alongside battery development, supported by CCGT (+30 GW) acting as backup, to host new 16 GW of DC capacity by 2050. France remains the most resilient location across all scenarios due to its robust nuclear baseline and the planned deployment of 9 GW of $H_2$CCGT. Similarly, the UK through a diversified portfolio of offshore wind (107 GW), nuclear (18 GW), and a flexible backup fleet of CCGT-CCS (10 GW) and BECCS (2 GW), mountains also relatively stable share in DCs growth investment.

The analysis of Poland reveals that the country initially experiences high growth by 2035, driven by relatively relaxed carbon constraints and a continued reliance on fossil fuels until 2045. This potential is further strengthened in the following decade through a more diversified energy mix, including the deployment of nuclear (17 GW) and $H_2$CCGT (3.7 GW). These factors transform Poland into a highly attractive destination for DC, capable of matching France's capacity in lift-off scenarios. In contrast, Italy presents a paradox. Despite possessing 45 GW of CCGT capacity by 2050, the country utilises natural gas as a baseload provider, with a 40% capacity factor by 2050, rather than as a flexible peaking resource. Consequently, despite solar and battery expansion, Italy remains gas-dependent, forcing the system to operate high-cost BECCS units to avoid emission violations. This high operational expenditure discourages DC investment after 2030. However, in pessimistic scenarios where DC capacity collapses between 2035 and 2045 and then recovers, Italy plays a crucial role in the subsequent recovery. In such cases, the system prefers to leverage existing gas infrastructure to meet additional DC demand rather than investing in new renewable or battery assets, then Italy can host +4 GW in more aggressive scenarios.

The most significant paradox occurs in the Nordic countries. These regions experience intense winter peaks and lower load factors compared to the rest of Europe, resulting in more volatile baseload profiles. While high availability of wind and hydro resources suggests a perpetual power surplus, these nations are highly susceptible during winter midday hours characterized by low wind. Due to strict carbon neutrality targets (nearly zero carbon post-2040), and lack of access to flexible, high-ramp sources such as gas forces these countries to rely heavily on electricity imports, a situation exacerbated by the presence of data centers. Although the region's abundant hydro (55 GW) and nuclear capacity (30 GW) play an important role in providing stable electricity, due to the electrification of other sector, as well as the ramp constraints pose significant barriers to DC expansion after 2030. This is particularly evident in the UK profile case, which is more stable and reduces the need for fast-ramp sources, where the DC profile requires a higher level of energy in the early morning hours and can be matched with regions that have high wind availability during that period (Supplementary Material Figure S.12 and S.13). Under this scenario, in both the baseline scenario and the lift-off scenarios, DC capacity could reach up to four times the 2030 capacity in the Nordic.

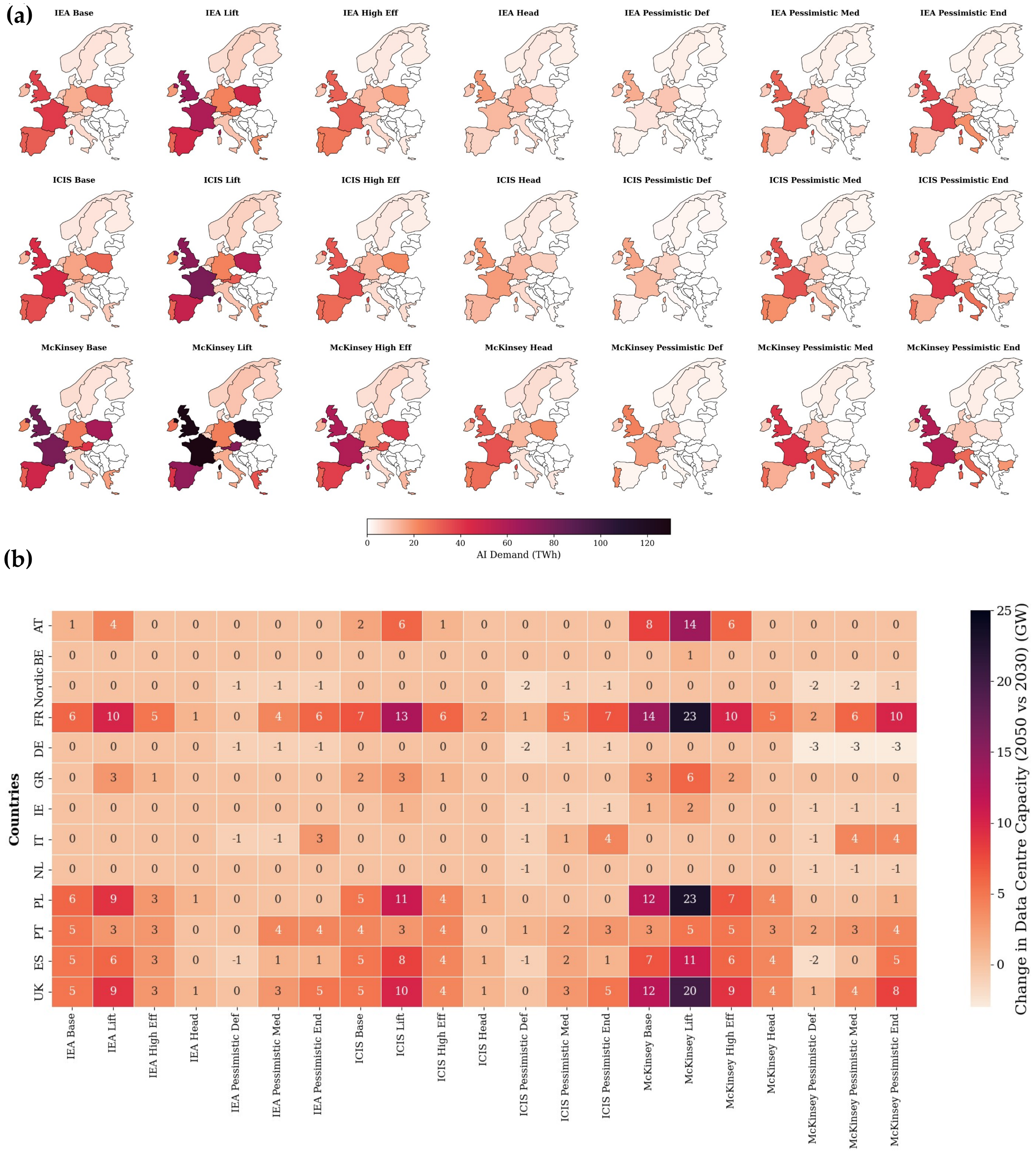

**Fig. 2. Geospatial distribution and capacity scaling of AI-driven DC demand across Europe by 2050. a,** Spatial mapping of annual AI electricity demand (TWh) across 21 scenarios, categorized by three primary demand trajectories (IEA, ICIS, and McKinsey). Each map represents the regional intensity of AI-related infrastructure under varying growth assumptions, ranging from "Pessimistic" to "Lift-off" scenarios. The colour scale indicates total demand, with dark purple regions identifying major hubs such as the UK, France, and Poland. **b,** Heatmap of projected AI-related installed capacity (GW) across 13 European countries/regions for each of the 21 scenarios. The grid displays the sensitivity of national infrastructure expansion to AI demand growth, with cell values representing gigawatts (GW).

**System response to stalled AI DC expansion**

For the sake of space, the subsequent analysis primarily focuses on the moderate ICIS trajectory. We conduct a counterfactual analysis to explore system responses if AI expansion stalls, comparing the ICIS Base scenario to three collapse variants (Def, Med, End), in which DC capacity contracts between 2035–2045. To further mimic real-life decision-making processes, the model adopts a myopic rolling horizon solution approach with a foresight of 5 year.

Fig. 3a shows that by 2050, peak load in the "Def" scenario is about 40 GWh lower than in ICIS Base. Despite this reduction, the share of firm generation is ~2.5% higher (y-axis). At the same time, the average Levelized Cost of Electricity (LCOE) in the collapse scenarios is €2–4 /MWh higher than in ICIS Base (bubble size). The drivers of this counterintuitive behaviour become clear in Fig. 3.c. All three collapse scenarios exhibit weaker deployment of renewables and storage, while the system retains roughly +40 GW of additional CCGT capacity compared with the ICIS base. Other firm technologies track the base trajectory in installed capacity, yet operate at higher utilisation rates, with capacity factors increasing by approximately: Nuclear: +4%, CCGT-CCS: +2%, $H_2$-CCGT: +3%, Biomass: +2%. Because substantial firm capacity was committed before the demand downturn (2035 in "Def", 2040 in "Med" and "End"), decommissioning these long-lifetime investments is financially suboptimal. Rather than expanding new renewables or storage, the system prioritises dispatching existing firm. This persistence increases system OPEX and drives the observed LCOE rise in Fig. 3a. A key observation is the increase in renewable spillage during collapse scenarios (Fig. 3b). Due to over-installed capacity in previous periods, the system cannot fully absorb renewable energy, leading to unavoidable curtailment. This highlights that, even with DC collapse, the system, particularly in pre-collapse firm capacity, cannot utilise renewable resources efficiently because of over-design.

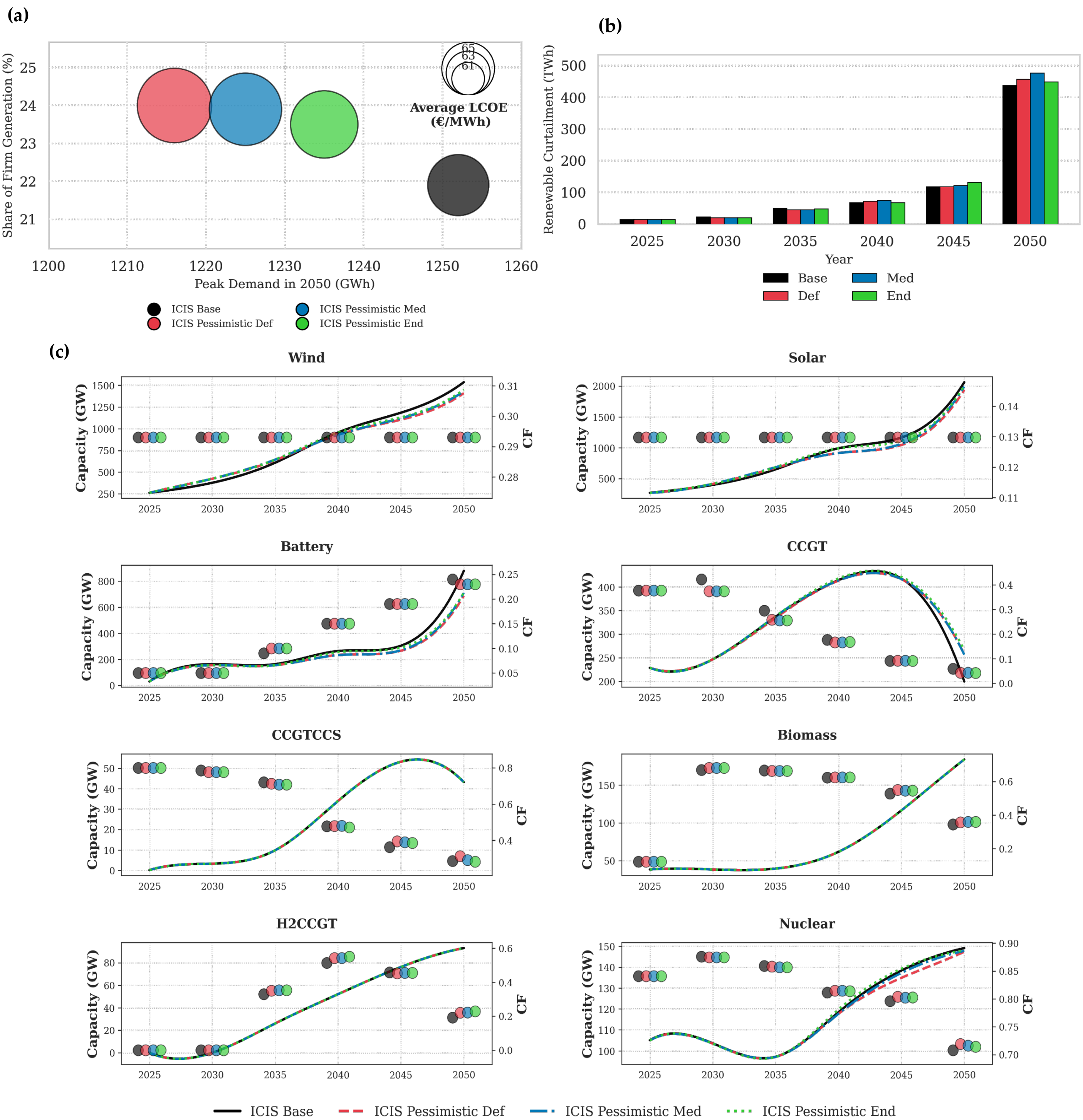

**Fig. 3. Impact of AI demand collapse on system adequacy, technology mix and (LCOE). a.** Trade-off between peak demand (GWh) and the share of firm generation (%) in 2050 for ICIS-based pessimistic scenarios (Deflation, Medium, and End) relative to the ICIS Base scenario. Bubble sizes represent the average system LCOE (€/MWh), illustrating a counterintuitive cost increase (€2–4/MWh) despite lower peak loads in collapse scenarios. **b.** Renewable curtailment trajectory. The panel highlights that the system fails to integrate renewable power during the collapse period, even with a reduction in installations compared to the base case. **c.** Temporal evolution (2025–2050) of installed capacity (GW, solid and dashed lines) and corresponding capacity factors (CF, represented by coloured circles) for eight power generation and storage technologies. The panels reveal a systemic shift in the pessimistic pathways: a reduction in variable renewable (Wind, Solar) and battery deployment is offset by the retention of approximately +40 GW of additional CCGT capacity and higher utilization rates of existing firm assets (Nuclear, Biomass, CCGT-CCS, and $H_2$-CCGT) to recover sunk costs from pre-collapse investments.

**Cross-border electricity dynamics and cost impacts under data-centre demand**

The integration of approximately 280 TWh of AI-driven demand by 2050 (under the moderate ICIS scenario) shows a fundamental shift in the European energy trade landscape. As demonstrated by Fig. 4.a, high-demand nations become net importers to sustain hyperscale operations. A primary example is the UK, where net imports increase from 17 TWh/yr in the Base scenario to 31 TWh/yr. This trade reversal is most pronounced in Spain, while it serves as a major exporter to France (over 50 TWh/yr) in the Base Case (no AI), the DC expansion transforms it into a net importer of 20 TWh/yr (Fig. 4.b). Denmark's export surplus shrinks from 35 TWh to about 10 TWh as rising domestic AI demand absorbs capacity, while France and Netherlands using nuclear baseload increase exports to Belgium and the UK. However, this systemic reliability comes at a significant environmental and economic cost.

The horizontal axis of Fig. 4.a reveals a sharp increase in cumulative emissions. Across Europe, the projected +77 $MtCO_{2eq}$ increase in emissions between 2030-2050, associated with approximately 53 GW of hyperscale DC capacity, is primarily driven by higher CCGT consumption (+300 TWh). In addition, an extra +8 TWh of coal generation in Central European countries, contributes significantly to rising emissions. Italy, hosting 1.6 GW of data-centre capacity by 2050, records an additional +10 Mt of $CO_2$, largely driven by +67 TWh of gas consumption. Germany and Poland emit around +5 $MtCO_{2eq}$ each, reflecting increased gas use (+36 TWh and +26 TWh, respectively) alongside additional coal generation (+4 TWh and +1 TWh). The Netherlands and Ireland also experience higher emissions (+5 $MtCO_{2eq}$ and +1 $MtCO_{2eq}$), associated with increased gas consumption of +20 TWh and +5 TWh, respectively. In the UK, emissions increase by approximately +7 $MtCO_{2eq}$, driven by +76 TWh of additional gas generation. However, the deployment of CCGTCCS, and BECCS, moderates the emissions intensity relative to capacity expansion compared to other countries. By contrast, France experiences a much smaller increase in emissions (below 1 $MtCO_{2eq}$), despite installing 8.5 GW of data-centre capacity. This is largely due to greater reliance on $H_2$-CCGT (+18 TWh) and nuclear generation (+120 TWh between 2030–2050).

The bubble radii in Fig. 4.a highlight the resulting economic costs. Meeting AI-driven demand by mobilising higher-cost resources, including coal, CCGT, nuclear, and $H_2$-CCGT, causes a noticeable rise in the LCOE. In the UK, LCOE increases from 125 to 160 €/MWh, while in Germany it increases from 106 to 120 €/MWh (Detailed of individual LCOE for all countries are reported in Supplementary Material Table S.6).

Fig. 4.b illustrates the changes in national generation and cross-border electricity exchanges after the introduction of AI demand, compared with the without AI (Country-level generation mixes and bilateral trade balances are reported in the Supplementary Material Figure S.10, S.11). Spain sees the largest generation increase (+100 TWh) from AI demand, mainly via solar, followed by France (+90 TWh), and UK and Poland (+60 TWh each), reflecting their roles as major AI load centres. AI demand also reshapes trade: Portugal shifts from exporting to importing ~10 TWh/yr from Spain, while UK imports ~5 TWh/yr from neighbouring systems, including Netherlands, Denmark, Norway, and Ireland. In the Nordics, AI demand shifts trade patterns: Sweden and Norway become modest importers (~2 TWh/yr from Germany), reversing their traditional export roles.

**(a)**

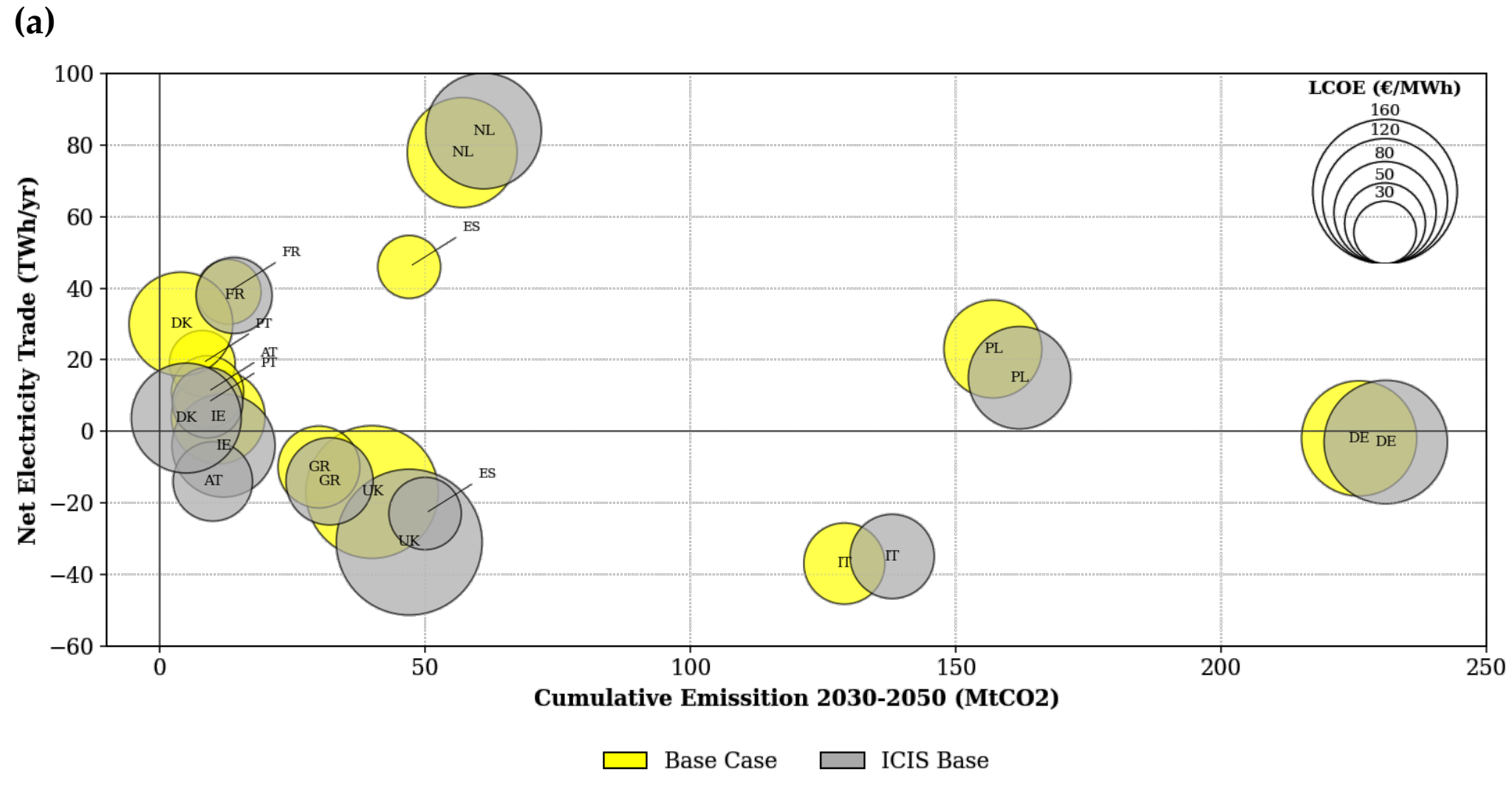

**(b)**

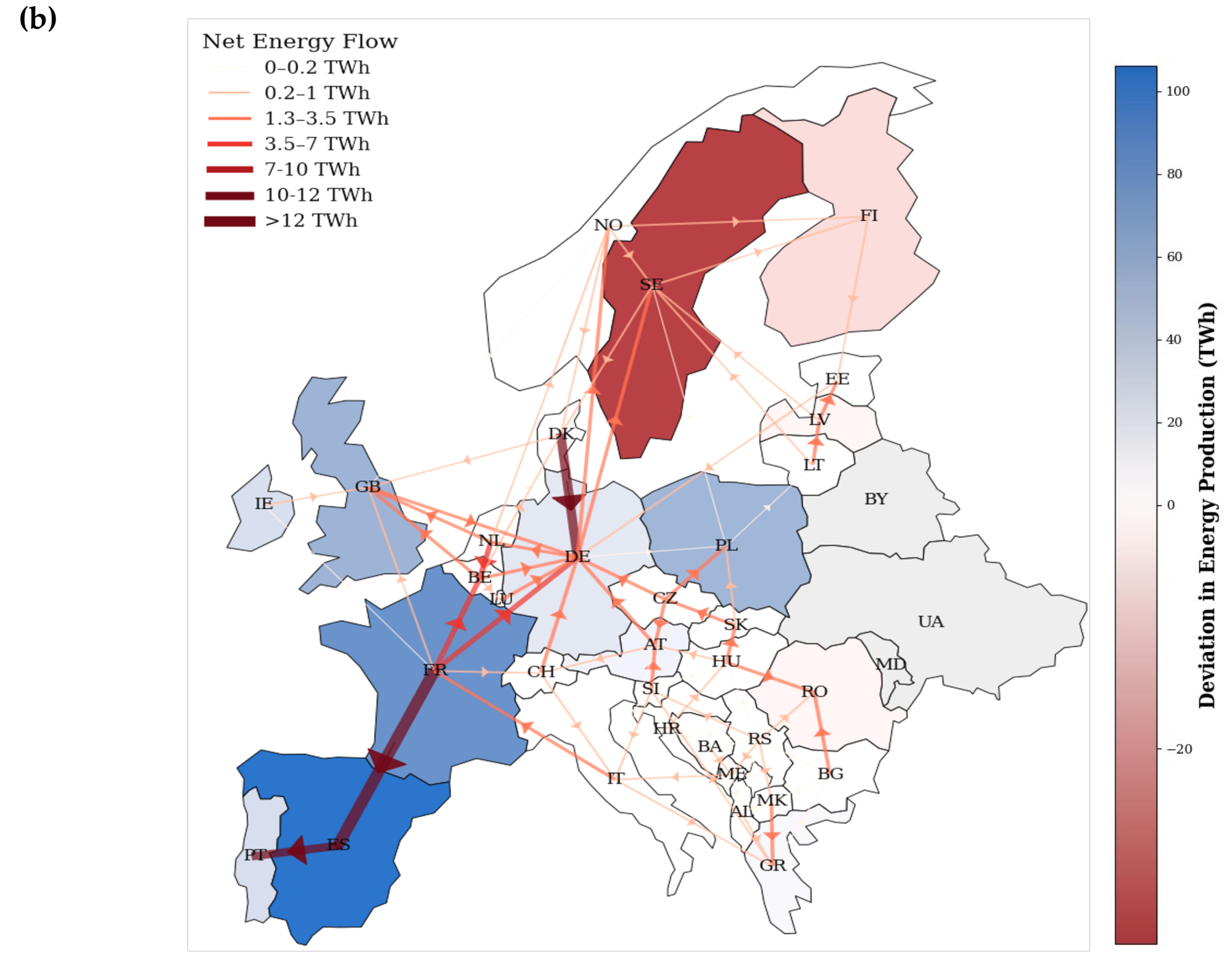

**Fig. 4 Impact of AI-driven demand on European electricity trade, carbon accumulation, and regional energy flows by 2050. a.** Multidimensional analysis of net electricity trade (TWh/yr) versus cumulative CO2 emissions (Mt CO2, 2030–2050) for 12 key European nations. Yellow bubbles represent the Base Case (without AI), while grey bubbles denote the moderate ICIS Base scenario. The y-axis identifies net exporters (positive) and net importers (negative), while bubble radii correspond to LCOE in each region. The shift from yellow to grey illustrates how AI-induced demand necessitates increased investment and thermal dispatch, subsequently shifting nations toward higher cumulative emissions and altered trade balances. **b.** Geospatial mapping of deviations in annual energy production (TWh, color-coded regions) and net energy flows (TWh, red arrows) in 2050. Blue shading indicates production increases (e.g., Spain, France), while red shading denotes production displacement (e.g., Sweden). Arrow thickness represents the magnitude of trade reversal or intensification.

**Influence of 24/7 AI demand on net demand profiles and dispatchable energy requirements**

This section examines the hourly impacts of DC operation. Fig. 5.a shows the load-duration curve (LDC) of net demand (demand minus non-dispatchable generation) for the baseline (Without AI) and ICIS scenarios under two DC load profiles, illustrating the number of hours the system relies on dispatchable resources. The US profile (52% load factor) exhibits higher variability, while the UK profile (67% load factor) is smoother and more uniform (Details are provided in Supplementary, S.1). While differences in the LDC over the entire year 2050 are visible across the three cases, the most significant effects appear in the focused section. In this area, AI operating under the US profile increase the number of hours requiring dispatchable resources by over 250 hours per year, whereas under the UK profile the increase exceeds 300 hours per year. This clearly demonstrates the impact of continuous, 24/7 AI demand on European electricity dispatch.

Fig. 5.b further illustrates the dispatchable energy requirements for two representative days, one in winter and one in summer, under the influence of DC for both load profiles. The US load profile exhibits greater temporal variability, with higher peaks and deeper off-peak troughs. Peak demand reaches approximately 1% above baseline in winter and 2% in summer, while off-peak demand falls to around 1% below baseline in both seasons. In contrast, the UK load profile is comparatively flat, showing less intra-day fluctuation. This stability leads the model to allocate a larger share of generation to nuclear capacity under the UK profile, with average nuclear dispatch of 331 MWh on a winter day and 120 MWh on a summer day.

The more stable demand structure in the UK case also reduces the need for flexible capacity. Consequently, the system installs less CCGT and $H_2$-CCGT capacity while investing more heavily in onshore wind, which can reliably supply a larger fraction of the relatively constant demand. Seasonal dispatch patterns further illustrate these structural differences. During winter, the model relies less on thermal balancing, with lower average outputs from CCGT (367 MWh) and $H_2$-CCGT (2050 MWh), as steady nuclear and wind cover a larger share of the demand. In summer, the pattern reverses: the UK profile requires greater flexible support, resulting in higher gas dispatch, including CCGT (1796 MWh) and $H_2$-CCGT (934 MWh).

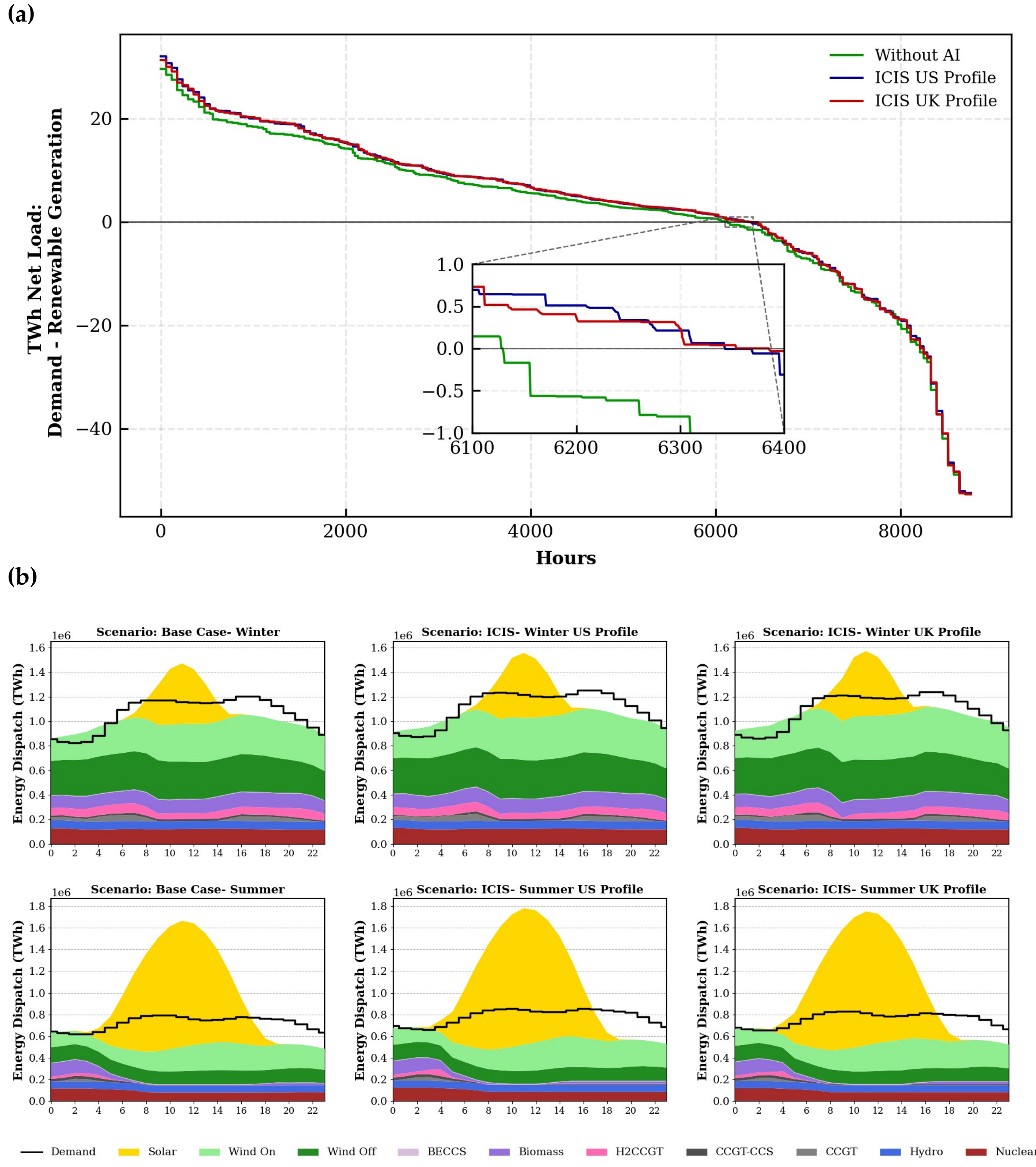

**Fig. 5 Temporal dispatch dynamics and impact of AI load profiles on system flexibility. a.** Net load-duration curves for the European power system in 2050, comparing the baseline (No-AI) with the ICIS scenario under two distinct AI demand profiles: U.S. (45% load factor) and UK (60% load factor). The zoomed-in section highlights the transition between surplus renewable generation and the requirement for dispatchable resources. AI-driven demand increases the annual duration of dispatchable resource necessity by 250 to 300 hours, depending on the load profile variability. **b.** Hourly energy dispatch stacks for representative winter (top) and summer (bottom) days across three cases: Base Case, ICIS-US, and ICIS-UK profiles.

**Hourly operational and market effects of DC demand across key European countries**

The spatial behaviour of the energy dispatch across countries highlights the impacts of DC more clearly on a representative winter day. In France, alongside the increase in nuclear generation, the presence of DC drives a noticeable midday rise in $H_2$-CCGT dispatch (+17 GWh) during a day (Fig. 6.a). In the UK, nuclear generation increases by an average of 1.1 GWh, biomass by 2.9 GWh on an average hourly basis across the day, and CCGT-CCS adds 6.4 GWh during the first six morning hours only where DC alone can add up to 5% to the hourly load during this period, leading to reduced electricity exports. In Germany, changes in $H_2$-CCGT dispatch around midday illustrate the compensatory role of DC in supporting system flexibility (+160 GWh daily). The system also requires additional imports in the early hours to meet demand (+57 GWh during 1:00-6:00). In Italy, the DC-driven demand significantly affects the load profile, forcing the system to generate up to 129 GWh more from CCGT in the early hours (1:00-9:00). In the Netherlands, changes in late-day load drive increased $H_2$-CCGT dispatch during this period (+40 GWh), while electricity exports decline (Detailed daily dispatch profiles for all countries, before and after AI integration, are provided in Supplementary Material Table S2 - S13).

AI deployment significantly alters 24-hour locational marginal prices (LMP). Countries with high AI penetration, such as France and the UK, exhibit higher prices in hours 16:00–22:00, while Germany, Netherlands, and Italy experience spikes in the early hours (Fig. 6.b). In contrast, countries without significant AI deployment experience lower electricity prices due to stagnant demand growth, expanding renewable generation, and the ability to export more power to neighbouring countries (negative prices).

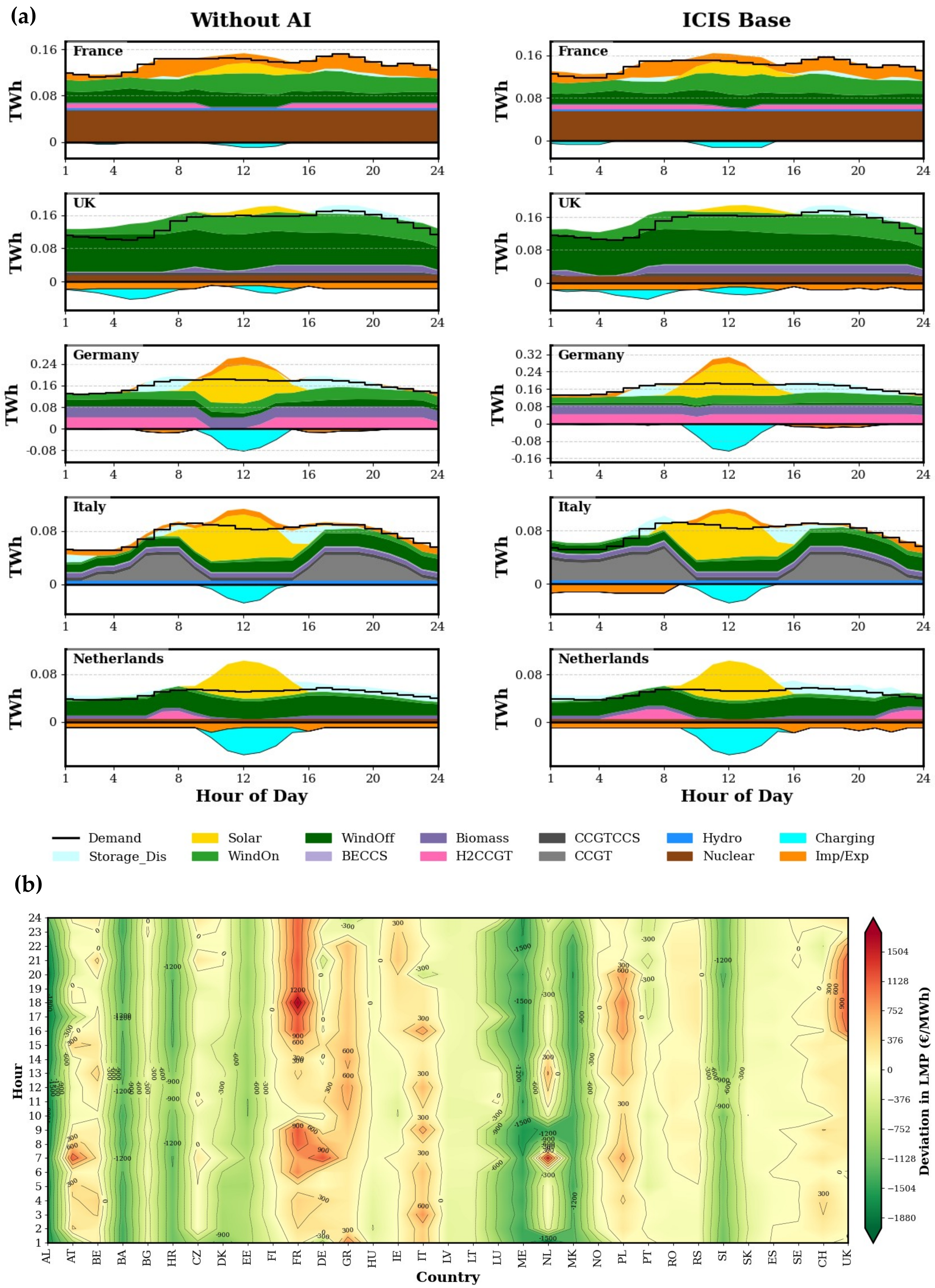

**Fig. 6. Hourly dispatch dynamics and marginal price volatility across key European hubs by 2050. a,** 24-hour electricity dispatch profiles for France (FR), United Kingdom (UK), Germany (DE), Italy (IT), and Netherlands

(NL) on a representative winter day, comparing the "Without AI" counterfactual (left) with the ICIS scenario (right). The stacks illustrate how AI-induced load increases trigger localized responses, such as elevated $H_2$-CCGT dispatch in France and Germany, and increased reliance on CCGT in Italy during early hours (1–9 a.m.). Dashed lines represent total demand, while shaded areas below the zero-line indicate energy storage charging or net exports. **b,** Spatiotemporal heatmap of 24-hour LMP deviations (€/MWh) following AI deployment. The colour gradient highlights sharp price escalations (red) in high-penetration regions like France and the UK during evening peaks (hours 16–24) and in Central Europe during off-peak hours. Conversely, price declines (green) are observed in non-AI regions, such as the Balkans, reflecting broader structural shifts in regional market equilibria.

**Impact of PUE on system outcomes: Sensitivity analysis**

This section examines how DC workloads and energy efficiency may affect the future European power system. The power usage effectiveness (PUE) is the key driver of electricity consumption in DCs. IEA projections indicate Europe's average PUE could fall from 1.45 today to 1.29 by 2030. Our Base Case assumes a conservative 1.25 post-2030, but uncertainties in technology or infrastructure upgrades require sensitivity analysis. Varying PUE from 1.4 to 1.1 impacts total cost, emissions, and grid peak demand (Fig. 7). If, contrary to expectations, PUE improves only to 1.4, representing a modest 5% reduction from current levels, total costs increase by approximately 8%, 7%, and 6% by 2030, 2040, and 2050, respectively, under the U.S. load profile. Conversely, achieving a PUE of 1.1 yields a 6%, and 4% cost reduction following the U.S., and UK profile, respectively by 2050.

A PUE of 1.4 raises carbon emissions by up to 8% in 2030 and 6% in 2040. The effect of efficiency varies with workload volatility: improving PUE to 1.1 under the U.S. profile reduces emissions by 5%, while the same improvement under the UK scenario cuts only 3%. This shows efficiency gains are roughly 1.6 times more effective for volatile AI workloads.

The 2050 peak-demand analysis shows that each 0.05 improvement in PUE lowers peak electricity demand by roughly 5 GWh for the U.S. profile, and 3 GWh for the UK profile. Without PUE improvements beyond 1.4, an additional 67 GW (U.S.) and 60 GW (UK) of generation capacity would be needed to maintain reliability. These results at the national level are even more pronounced potentially reducing the need for up to 50 GW of new capacity by 2050 (Supplementary Material Figure S.12 and S.13).

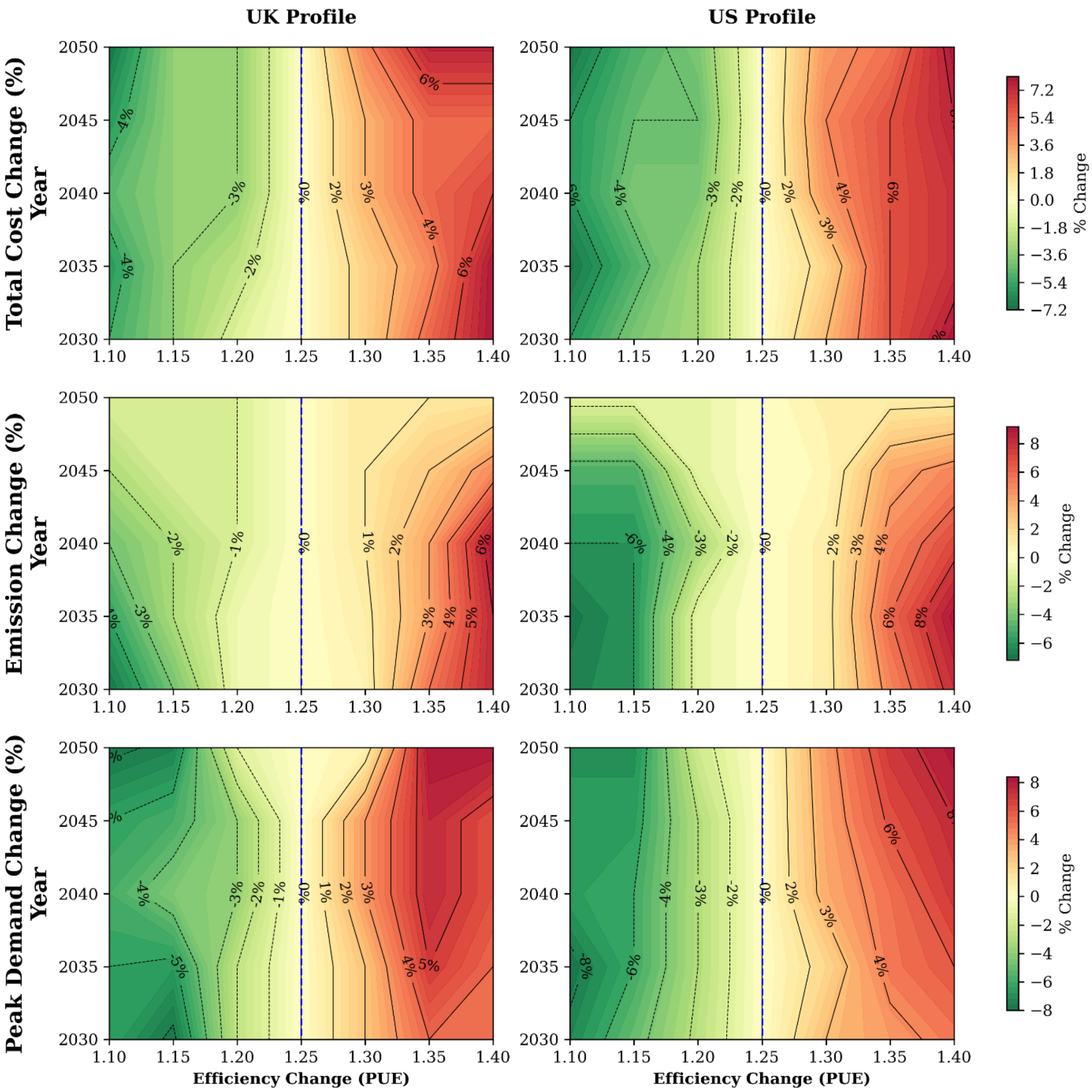

**Fig. 7. Multi-dimensional sensitivity analysis of AI workload efficiency on grid metrics.** The contour plots illustrate the percentage change in total system cost (top), cumulative emissions (middle), and peak demand (bottom) as a function of PUE (x-axis) and temporal progression from 2030 to 2050 (y-axis). The left column represents the UK steady profile, while the right column depicts the US volatile profile. The black isolines represent constant percentage intervals, where the slope and density of these lines indicate the sensitivity of the metric to efficiency changes over time. The vertical blue dashed line marks the 1.25 PUE threshold, separating the beneficial efficiency zone (green, PUE < 1.25) from the systemic penalty zone (red, PUE > 1.25).

## Discussion

This study shows that rapid DC expansion alters the structure of Europe's power system transition. The primary constraint is not only annual energy supply, but firm capacity adequacy. Under current infrastructure trajectories, this structural shift induces a temporary reversal in fossil phase-out dynamics and amplifies transition risks during the 2030s.

Across 21 growth scenarios, DCs expansion systematically increases reliance on firm thermal resources. In the medium term, gas emerges as the dominant stabilizing technology, while coal and oil persist beyond their expected phase-out window. In lift-off scenarios, coal remains in the generation mix until 2045 rather than exiting by 2035 as in non-AI baselines. Seventeen scenarios exhibit emission overshoot between 2030 and 2045. AI-driven power demand does not prevent the achievement of long-term net-zero targets. However, it substantially increases the likelihood of exceeding Europe's interim emission targets during the 2030s. This deviation reflects a temporary gap between the pace of DC demand expansion and the deployment of firm, low-carbon dispatchable capacity required to maintain system adequacy.

Our results indicate that Europe's baseline infrastructure cannot accommodate AI growth without structural expansion, requiring an additional 73–723 GW of capacity across pessimistic to lift-off scenarios, of which roughly 70% consists of renewables and grid-level storage and about 30% of firm dispatchable resources. Under existing capacity assumptions, the system must either expand firm resources, primarily CCGT, $H_2$CCGT, and nuclear, or accept a delayed fossil phase-out. Gas shifts from a marginal balancing role to a structural adequacy anchor during the 2030s. As carbon constraints tighten in the 2040s, its function becomes increasingly backup-oriented, yet early capacity additions generate inertia and potential lock-in.

The geographic allocation of DC capacity further underscores that annual renewable abundance is not a sufficient siting criterion. Regions commonly perceived as optimal due to surplus renewable generation, particularly Nordic systems, face structural adequacy constraints. Despite high yearly energy balances, limited dispatchable capacity exposes these systems to shortfall risk during critical hours. The binding constraint for hyperscale DC integration is therefore not energy volume but firm, low-carbon dispatchability. Hence, transforming them into future DC hubs will strongly depend on flexible generation. Conversely, regions with diversified portfolios combining renewables, nuclear baseload, and flexible gas assets demonstrate greater resilience under worst-day conditions. In this context, $H_2$CCGT technologies play a pivotal bridging role in maintaining system stability while preserving long-term decarbonisation optionality.

Our results also identify a distinct risk in the case of DCs demand growth collapse. In scenarios where AI demand weakens after 2035, early investments in firm capacity result in capacity deadlocks. Oversized thermal fleets elevate system-level LCOE and prolong fossil backup utilisation even under slower demand growth. This creates a new category of transition risk: digital overbuild lock-in, whereby anticipatory infrastructure expansion limits future flexibility.

Improving efficiency and enhancing workload flexibility could become two key parameters in the future of DCs, repositioning them not as an additional burden on the grid but as a source of flexibility. In doing so, DCs could support renewable integration, reduce peak demand, reduce total installed capacity and lower overall emissions.

We demonstrated that the temporal behaviour of DC workloads strongly influences the energy dispatch structure and, consequently, the system's flexibility requirements, emissions trajectory, and dependence on different generation technologies: A more constant load profile increases the system's reliance on nuclear and natural gas generation, while also encouraging greater wind capacity expansion and lower solar deployment. At the same time, the reduced need for flexible balancing resources leads to a lower role for $H_2$-CCGT generation in the electricity mix.

Several limitations warrant further investigation. Land-use constraints for hyperscale facilities, local distribution bottlenecks, and connection queues are not explicitly modelled and could materially affect feasible deployment rates. While sensitivity analyses account for improvements in PUE, rapid advances in hardware efficiency and cooling technologies may alter future load intensity and siting patterns. Finally, electric vehicles (EVs) and demand response program dynamics as flexible options are not incorporated in our model. EVs are represented as stationary storage rather than dynamic vehicle-to-grid actors, potentially underestimating system flexibility. Additionally, in this study, demand response programs, both at the network level and from the perspective of DCs, which could contribute to improving network flexibility, have not been examined and are left for future work. Also, it is important to note that, in this study, AI has been considered only from the perspective of electricity consumption and the additional pressure it may impose on the power system. However, it is well recognized that AI has significant and widespread applications in improving efficiency and reducing electricity consumption across various sectors, such as forecasting methods, cost reduction, and numerous economic optimisations.

Overall, even in the most optimistic scenarios where AI-driven data centers account for around 9% of Europe's total electricity consumption by 2050, the challenge is not merely one of aggregate energy demand; rather, it reshapes the adequacy architecture of Europe's energy transition. Managing this shift requires aligning digital expansion with accelerated deployment of firm low-carbon capacity to avoid temporary fossil reliance and long-lived structural lock-in.

## Methods

### Problem description

A full description of the problem, including all assumptions, formulations, and input data, is provided in the Supplementary Material. Here, we present a brief overview of the problem, highlighting the key features of the optimisation model and its main assumptions, to facilitate the interpretation of the results presented in the preceding sections.

#### Long-term European power grid expansion model

The long-term generation and transmission expansion planning problem for the integrated European power system is formulated as a multi-period, spatially explicit MILP model that simultaneously optimises investment and operational decisions. The objective function minimises total system cost, comprising both capital and operational expenditures. Capital costs include investments in new generation capacity, energy storage systems, and cross-border transmission infrastructure. Operational costs account for variable generation and storage costs, fuel consumption, fixed operation and maintenance costs, carbon cost, penalties for loss of load (value of lost load, VOLL), and penalties for renewable energy curtailment. The full details of problem formulation including objective function and constraints are provided in Supplementary Material, Section S. 3.

The spatial resolution of the model covers 33 European countries, including EU Member States as well as the UK, Norway, Albania, Bosnia and Herzegovina, Serbia, North Macedonia, and Switzerland, with real existing cross-border transmission interconnections explicitly represented [36]. A total of 177 real transmission projects between countries are considered as candidate cross-border capacity expansions that can be invested in over the planning horizon [37]. To reduce the computational complexity of the large-scale optimisation problem, the planning horizon to 2050 is divided into six multi-year periods with five-year strategic time steps. At the annual level, electricity demand and the availability of onshore wind, offshore wind, and solar resources were derived for the 33 countries under study, resulting in a dataset comprising **33 × 4 × 8,760** hourly time series. These annual profiles for each planning years, reduced to 12 representative days with hourly resolution (24 h) using k-medoids clustering [38,39],with three representative days selected per season to adequately capture seasonal variability (Supplementary Material, Figure S.3). Then, operational decisions, such as generation output, storage charging and discharging, and cross-border power flows, are determined at an hourly resolution for the representative days, with weighting factors used to scale results to annual system performance.

**Model input & assumptions**

The proposed optimisation framework is designed to employ state-of-the-art methodologies and the most recent datasets to enhance the realism of the results. In this study, the Global Ambition scenario from the TYNDP developed by the ENTSO-E was selected as the reference framework for the analysis [40]. The Global Ambitious scenario assumes strong and coordinated global and European climate actions aimed at achieving deep decarbonisation of the energy system. It reflects a pathway consistent with the European's objective of reaching climate neutrality (net-zero greenhouse gas emissions) by 2050. This scenario is characterised by the rapid expansion of renewable energy sources, widespread electrification across transport, industry, and heating sectors, as well as the large-scale deployment of low-carbon technologies such as green hydrogen and energy storage. Given that the primary objective of this research is to assess DC power demand expansion under conditions compatible with the EU's net-zero target by 2050, the Global Ambitious scenario provides an appropriate and policy-consistent analytical framework. While the core structure and key inputs, such as technology build rates at the country level, carbon reduction policies, fuel and carbon prices, are aligned with the TYNDP methodology, the model departs from TYNDP by jointly solving capacity expansion and system dispatch, thereby requiring the explicit representation of multiple technical constraints to ensure system adequacy. Unlike TYNDP scenarios, where decarbonisation is largely enforced through minimum renewable capacity constraints and thermal capacity expansion (including gas) is effectively prohibited after 2030 in Global Ambition scenarios, our framework allows endogenous investment decisions subject to adequacy requirements. Recent adequacy assessments based on inputs from European TSOs indicate that several countries, including the UK, Finland, Portugal, Slovenia, Italy, and Poland, retain potential for gas capacity additions between 2030 and 2035, with expectations that such resources may remain necessary in subsequent decades to maintain system reliability [41]. Accordingly, the proposed model does not impose minimum build-rate constraints for renewable technologies, and endogenous development of CCGT capacity is permitted.

Moreover, while TYNDP does not explicitly enforce must-run constraints for thermal units beyond 2030, neglecting such constraints, particularly minimum generation and ramping limits, can lead to unrealistic dispatch outcomes, especially for technologies such as nuclear power. Although unit commitment is not modelled in our model, must-run constraints are incorporated for nuclear, resulting in higher nuclear generation compared with TYNDP projections (by approximately 150 TWh), reflecting more realistic operational behaviour.

One of the key differences between the proposed model and previous approaches is that system adequacy is explicitly supported through a spinning reserve constraint, parameterized using technology-specific derating factors, as formulated in Eq. (1).

$$\sum_j \sum_c de_j \times CAP_{j,c,t} \geq (1+\gamma) \times PeakDemand_t \qquad \forall\, t \qquad (1)$$

Where, $de_j$ denotes the derating factor, and $\gamma$ represents the marginal reserve, which is set to 8 % in the model, $CAP_{j,c,t}$ is the total capacity of technology $j$ in country $c$, at year $t$, and $PeakDemand_t$ is the whole-system peak demand at year $t$. It should be noted that the system reliability is ensured using a loss of load expectation (LOLE) constraint, limited to three hours per year for each country [41].

In addition, whereas many European power system expansion models include penalties VOLL, the cost of renewable curtailment is often neglected. Given projections indicating that renewable curtailment could exceed 310 TWh by 2040 due to grid congestion [42] (outside the scope of this study), a curtailment penalty of 200 €/MWh is explicitly included in the objective function. To further enhance realism, technology-specific land availability constraints based on empirical data are applied for renewable deployment [43,44].

Beyond the primary assumptions, the model accounts for country-specific distribution loss [41], and cross-border transmission losses, as well as technology-specific self-consumption / self-charge for battery storage, to provide a more realistic representation of system operation (S. 3 in Supplementary Materials).

## Data centre capacity scenario pathway

Although numerous quantitative and qualitative forecasts regarding the explosive growth of DCs, particularly AI-workload-driven DCs, have been reported across various studies and industry reports in the scope of Europe, the fundamental question of how this growth will be distributed among countries has not been sufficiently examined, either from a predictive or an optimal allocation perspective. In other words, the extent to which different European regions can host hyperscale DC capacity, and how such capacities may affect Europe's long-term energy transition strategies, remains an open and critical question. Addressing this issue requires a detailed analysis of national and supranational AI-related policies at the regional level.

In this study, using 21 pathway scenarios derived from three major databases: IEA, McKinsey, and ICIS, we estimate the total capacity of AI-based DCs across Europe up to 2050. The IEA, ICIS, and McKinsey databases provide different annual growth projections for total DC capacity across Europe. Nevertheless, all forecasts consistently indicate that approximately 50–70% of future DC capacity will be dedicated to AI workloads [10,45], encompassing all hyperscale DCs as well as a substantial share of colocation facilities with capacities exceeding 10 MW [46]. According to the IEA, total DC capacity in Europe is projected to increase from the current 16 GW to 27 GW by 2030, corresponding to a CAGR of 11%. ICIS projects a higher annual growth rate of 15%, resulting in a total capacity of 30.4 GW by 2030, while McKinsey forecasts an even more aggressive growth scenario, reaching 35 GW by 2030 with an annual growth rate of 20% (However, as discussed in the preceding sections, ICIS projections for energy demand generally exceed those of McKinsey, reflecting a more load factor for DC). Beyond its base scenario, the IEA provides three additional global scenarios for DC capacity up to 2030: Lift-Off, High-Efficiency, and Headwinds, corresponding to annual growth rates of approximately 21%, 6%, and 3%, respectively [15,47]. Extending these scenarios to the projections provided by ICIS and McKinsey yields a broader range of possible outcomes that must be considered in analytical studies. Given regional heterogeneity and substantial political and geopolitical uncertainties, both more explosive growth trajectories and potential market collapses must be accounted for. In addition to these scenarios, and to enhance the robustness of the analysis, we assume that existing trends over 5- or 10-year intervals may be disrupted by various shocks. Accordingly, three pessimistic scenarios, Deflection, Mid, and End, are introduced. A major limitation across existing studies is that most projections are restricted to 2030, with the IEA explicitly noting that forecasts beyond this horizon entail significant uncertainty. Nevertheless, by adopting the proposed growth patterns (rapid growth, moderate growth, stagnation, and collapse across different time intervals), it can be argued that nearly all plausible future market conditions are captured within this framework. Consequently, we extend these scenario pathways to 2050 and examine the overall trajectory of the European DC market under 21 distinct scenarios. For example, in the IEA Base scenario, annual growth is assumed to be 11% up to 2030, 10% for the period 2030–2035, 5% from 2035 to 2040, and 2.5% from 2040 to 2050. Annual capacity targets aggregated over five-year intervals are used as European targets in optimisation model. Exact numerical values and the methodology used to calculate the scenario pathways are fully described in the Supplementary Material, Section S.2.

**Regional data centre capacity**

Although obtaining precise current and future capacity data at the country level is highly challenging and no standardized forecasting methodology exists for installed DC capacity, the EUDCA provides current capacity estimates and annual projections for hyperscale and colocation DCs across Europe according to colocation and hyperscale DC database, Pb7 Research[48]. The European DC database is compiled and updated by Pb7 Research using a combination of existing sources, online searches, and data verification procedures. Where direct data are unavailable, estimates are derived from alternative information, such as rack counts or building layouts. These estimates are derived from data collection and interviews conducted with owners and operators of DCs across Europe between November and December 2024. Overall, the reported capacities show reasonable consistency with the figures presented in the IEA report. It should be noted that all reported values refer to IT capacity. The exact IT capacity for hyperscale and colocation DC for all European countries in 2024 is provide in the Supplementary Material, Table S.1. The credibility of these figures can be validated by comparison with the limited number of country-specific reports available for markets such as the UK, Germany, and the Netherlands, Belgium and Poland [49-52].

## Hyperscale data centre load profile

In addition to installed DC capacity, both the associated energy consumption and the contribution of DCs to each country's total electricity demand must be taken into account [53,54]. To enhance model robustness and validation, two distinct types of hyperscale data-centre load profiles are incorporated. The first is derived from half-hourly utilisation data for 87 DCs in the United Kingdom covering the period 2023 to July 2025[55]. By clustering these load profiles, all DCs are classified into three voltage levels: low, high, and extra-high voltage. Since hyperscale DCs with capacities exceeding 10 MW (at minimum) require high-voltage connections (110, 230, or 400 kV) and are typically connected to the transmission or sub-transmission network [46,56], we derive a representative load profile from 57 DCs connected at high-voltage levels using a clustering-based selection method.

The second consists of high-resolution load profiles from hyperscale DCs located in Texas [57]. The historical hourly load profiles for both datasets are presented in the Supplementary Material, Figure S.2. A key distinction between the two lies in their temporal variability: the UK-based DCs exhibit a relatively stable load profile with a load factor consistently exceeding 67%. In contrast, the Texas hyperscale data-centre profiles display more pronounced daily and seasonal fluctuations throughout the year, with a load factor consistently exceeding 52%.

Hence, the energy consumption by DCs is calculated as follows:

$$DemandAI_{t,c,k,h} = DCProfile_{k,h} \times CapAI_{t,c} \times PUE_t \qquad \forall\, t, c, k, h \qquad (2)$$

Where $DemandAI_{t,c,\mathrm{k},h}$ is a variable representing DC electricity demand in year $t$, country $c$, cluster $k$ and hour $h$, $CapAI_{t,c}$ denotes the optimal DC IT capacity in year $t$ and country $c$, is the PUE at year $t$ followed by the IEA projections (1.45 in 2024 and 1.29 in 2030) [47], and $DCProfile_{k,h}$ is the normalised hourly demand profile used as model input in cluster $k$ and hour $h$.